\documentclass[final]{siamltex}

\usepackage[colorlinks=true,urlcolor=blue]{hyperref}
\usepackage{amstext}
\usepackage{amssymb}
\usepackage{amsmath}
\usepackage{mathrsfs}
\usepackage{array}
\usepackage{enumerate}
\usepackage{verbatim}
\usepackage{url}

\usepackage[]{graphicx}               
\usepackage{color}
\usepackage{appendix}

\newcommand{\reals}{\mathbb{R}}
\newcommand{\DD}{\displaystyle}

\title{Electrostatic forces on charged surfaces of bilayer lipid membranes \thanks{This 
        work was partially supported by the Simons Foundation.}} 

\author{Michael Mikucki\thanks{Department of Mathematics,
        Colorado State University, Fort Collins, Colorado, 80523-1874
        ({\tt mikucki@math.colostate.edu}).}
        \and Y.~C.~Zhou\thanks{Department of Mathematics,
        Colorado State University, Fort Collins, Colorado, 80523-1874
        ({\tt yzhou@math.colostate.edu}).}}

\begin{document}

\maketitle

\begin{abstract}
Simulating protein-membrane interactions is an important and dynamic area of research. A proper definition of electrostatic forces
on membrane surfaces is necessary for developing electromechanical models of protein-membrane interactions. Here, we model the bilayer 
membrane as a continuum with general continuous distributions of lipids charges on membrane surfaces. A new electrostatic potential 
energy functional is then defined for this solvated protein-membrane system. Key geometrical transformation properties 
of the membrane surfaces under a smooth velocity field allow us to apply the Hadamard-Zol\'{e}sio structure theorem,
and the electrostatic forces on membrane surfaces can be computed as the shape derivative of the electrostatic energy functional.
\end{abstract}



\pagestyle{myheadings}
\thispagestyle{plain}
\markboth{MICHAEL MIKUCKI AND Y. C. ZHOU}{ELECTROSTATIC FORCE ON CHARGED DIELECTRIC INTERFACES}

\section{Introduction}
This paper concerns the mathematically rigorous and physically justifiable definition
of electrostatic forces on the surfaces of lipid bilayer membranes within the framework
of implicit solvent and continuum models of charged lipids. Lipid bilayer membranes serve 
as the boundaries of cells and many cell organelles and control the exchange of ions, 
nutrient particles, and metabolic products between the enclosed structures and the surrounding aqueous 
environment. Bilayer membranes function by stretching, bending, merging and separating to control gate specific 
channels or to wrap/unwrap particles. These deformations are precisely regulated by various proteins and 
other macromolecules, each with its own specific function. Among all types of intermolecular interactions that 
drive membrane deformations, the electrostatic interaction is ubiquitous, since proteins and lipid bilayer 
membranes are always charged under physiological conditions. The electrostatic interaction is especially important
if the membrane dynamics involves the lateral diffusion of charged lipids and electrostatic association of 
proteins. Computation of the electrostatic forces on the membrane is therefore of critical importance 
for the quantitative study of membrane dynamics and related cellular activities. Interested readers are referred 
to \cite{ChoW2005review,HonigB1986a,MclaughlinS2005a,Murray2001,WhiteS1999a} for more thorough 
discussions of the membrane dynamics and protein-membrane electrostatic interactions.

Computation of electrostatic forces on bilayer membranes by summing up pair-wise Coulombic interactions 
requires the full molecular details of the system and involves all particles in the computation; therefore, it is 
challenging to use this technique to study membrane dynamics at biologically relevant spatial and time scales. 
Implicit solvent models have been introduced so that the averaged behavior of highly dynamic solvent molecules can be described 
as a structure-less continuum \cite{Feig2004}. This simplification greatly reduces the degrees of freedom in simulations.  Such a continuum membrane model is 
necessary for simulating membrane dynamics over a region larger than 100\AA~in space and longer than a microsecond 
in time \cite{Choe2008,KhelashviliG2009a,Tang2008,Zhou2010}. 
A rigorous definition of the electrostatic force on biomolecules within the framework of implicit solvent models
has been a focus of computational biologists and applied 
mathematicians for decades \cite{Gilson93,Li2011,Sharp1990a,Zhou2008}. On one hand, the classical electromagnetic theory shows that
the electrostatic force on a dielectric interface can be related to the jump of the Maxwell stress tensor (MST) across the interface \cite{Stratton1941}. On the other hand, biophysicists and applied mathematicians usually compute the electrostatic force as the 
variation of some free energy functional with respect to the position of the dielectric interface. Because the same 
equilibrium state of a system can be given by the extremization of different energy functionals (e.g. \cite{MaggsA2002a,MaggsA2012a}), and 
variations of these energy functional (hence the simulated dynamics of the system) may not be the same, caution shall be
taken in defining an energy functional that can give correct dynamics and correct equivalent state. The equivalence of 
the force presentations obtained from MST and from enery minimization can be conveniently used for this purpose. 
In this work, we will verify that the electrostatic forces derived from our electrostatic energy functional matches
exactly with that defined by the Maxwell stress tensor.

Recently, Li \textit{et al.} obtained an elegant derivation of 
the electrostatic force on molecular surfaces by computing the shape derivative of the electrostatic potential 
energy \cite{Li2011}. His work assumes that the dielectric interface is uncharged, and therefore it is not applicable to bilayer 
membranes if they are modeled with continuous distribution of charged lipids on the surfaces. 
The variation of charge densities on dielectric interfaces must be considered in deriving the electrostatic force, and the 
model in \cite{Li2011} must be modified before it can be applied to protein-membrane interactions. Here in this work, we adopt the shape
derivative approach for modeling bilayer membranes with variable surface charge densities. Under a special condition, this surface charge density follows the Boltzmann distribution. The deformation of the bilayer membrane is governed by a smooth velocity field that vanishes at distances away from the membrane surfaces.  The surface transformation is approximated by a Taylor expansion of the velocity field.  
We find that the time derivative of the Jacobian of this surface transformation equals the
surface divergence of the velocity field on membrane surfaces (see a similar result in Lemma 4.2 in \cite{Simon1980}). Thanks to this 
essential geometric property, 
the Hadamard-Zol\'{e}sio structure theorem \cite{DelfourM1992a} holds true for our model, and thus the shape derivative approach of computing the electrostatic 
force can be applied. Our main result, the dielectric boundary force, is extracted from the shape derivative of the electrostatic 
potential energy. Finally, we justify the physical description of the force by expressing it in terms of the Maxwell stress tensor (MST).
Shape differentiation with respect to the domain or the boundary has been studied for various optimization 
problems \cite{ButtazzoG1991a,AllaireG2001a,NovruziA2002a,EpplerK2007a,CanelasA2009a},
and a general functional framework exists \cite{Simon1980}. This framework might be extended to our problem by matching 
the conditions across the internal dielectric interface. 

The major limitation of our model, as with all sharp-interface models, is that it does not describe topological changes to the membrane such as those occur during the processes of membrane merging and separating. The phase-field method is capable of modeling these phenomenon. We refer readers to \cite{Du2004, Du2006, Teigen2011,GavishN2011a,DaiS2013a} 
for applications of this method.  Furthermore, as the dielectric boundary force described by this model is induced by proteins that are proximal to but not in direct contact with the membrane surface, these proteins must be at a minimum distance $d>0$ away from the membrane to preventing the protein atoms from overlapping with the membrane surface.  

\begin{figure}[t] 
\centering
\includegraphics[height=6cm]{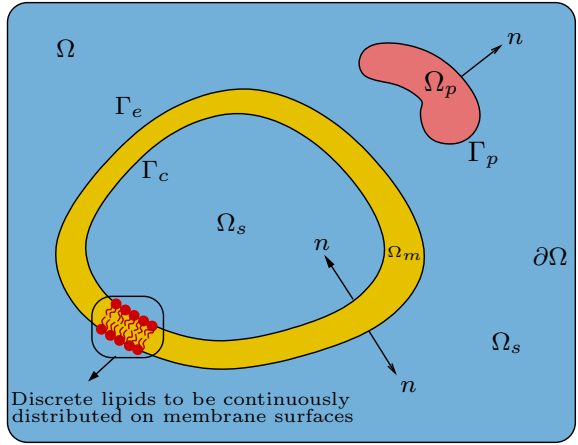}
\caption{Mathematical description of protein-membrane system.  The containment domain is denoted by $\Omega$.  
The volume of the lipid bilayer is $\Omega_m$, and the exoplasmic (exterior) and cytosolic (interior) faces of 
the membrane are $\Gamma_e$ and $\Gamma_c$, respectively. The volume enclosed by the membrane and the aqueous surrounding environment 
are both denoted $\Omega_s$.  The protein is $\Omega_p$ with surface $\Gamma_p$.  The unit outward normal to 
any surface $\Gamma$ is $n$, always pointing toward the solvent $\Omega_s$.  Lipids are drawn in the bottom left corner for illustration but their atomic 
details are neglected in the model.  Note that the protein may be located in the solvent region inside $\Gamma_c$ in some cases.} 
\label{fig:system} 
\end{figure}

We consider the following model problem for the protein-membrane interaction.
Define $\Omega \subset \reals^3$ as the region 
containing the entire membrane-protein system.  Let $\Omega_p$ denote the volume of the protein and $\Omega_m$ be the volume of the membrane.  We emphasize that $\Omega_m \subset \reals^3$ is a volume with finite thickness rather than a two dimensional manifold without thickness.
The solvent region is $\Omega_s$ which includes both the external surrounding environment to the bilayer membrane and the volume enclosed by the membrane.  The 
boundary separating the protein and the solvent is the manifold $\Gamma_p \subset \reals^2$.  
The membrane has two boundaries to the solvent.  The interior boundary is known as the cytosolic face and is denoted $\Gamma_c \subset \reals^2$.  
The exterior boundary is known as the exoplasmic face and is denoted $\Gamma_e \subset \reals^2$.  The exterior boundary of the containment 
domain is $\partial \Omega$.  Each boundary $\Gamma_p$, $\Gamma_c$, $\Gamma_e$, and $\partial \Omega$ is assumed to be smooth.  The atomic details of the lipids that compose the membrane are neglected for the continuum model.
The unit outward normal to any boundary $\Gamma$ is denoted as $n$.  A cross section of the protein-membrane model is illustrated in 
Figure \ref{fig:system}.  

The total potential energy of bilayer membrane is expressed in two components.  First, the classical mechanical bending energy depends only on the 
curvatures of the bilayer membrane.  Following the work of Canham \cite{Canham1970}, Helfrich \cite{Helfrich1973} and Evans \cite{Evans1974},  
for example, this bending energy is given by
\begin{equation}\label{eq:bending_energy}
E[\Gamma] = \int_{\Gamma} \left(\frac{1}{2} \mathscr{K}_C (2H-C_0)^2 +\mathscr{K}_G K \right)  dS, 
\end{equation}
where $\mathscr{K}_C$ and $\mathscr{K}_G$ are the bending modulus and Gaussian saddle-splay modulus, respectively, $H$ is the mean curvature, $C_0$ is the spontaneous curvature, and $K$ is the Gaussian curvature \cite{Feng2006}.  In addition to the bending energy, the membrane is under an external potential force induced by the proteins. The potential energy from this 
interaction is added to $E[\Gamma]$ to get the total potential energy of the system, $\Pi[\Gamma]$,
\begin{equation}\label{eq:total_energy}
\Pi[\Gamma] = E[\Gamma] + G[\Gamma],
\end{equation}
where $G[\Gamma]$ is the electrostatic potential energy from the protein-membrane interaction which will be defined in \eqref{eq:G_gen} below. The equilibrium position of the bilayer membrane is determined by the membrane surfaces which extremize the bending energy and the electrostatic potential energy.  
The extremization of \eqref{eq:total_energy} gives rise to the bending curvature equation of $\Gamma$, 
\begin{equation}
\delta_\Gamma \Pi[\Gamma] = \delta_\Gamma E[\Gamma] + \delta_\Gamma G[\Gamma].
\end{equation}
The variation of the classical bending energy has been calculated in \cite{Feng2006} to be 
\begin{equation}
\delta_\Gamma E[\Gamma] = \int_\Gamma \left[ \mathscr{K}_C 2(2H-C_0)^2 \delta H \sqrt{a} + \mathscr{K}_C\frac{1}{2}(2H-C_0)^2\delta\sqrt{a} \right] (1/\sqrt{a})\;dS.
\end{equation}
where $a$ is the determinant of the covariant metric tensor, and $s_1$ and $s_2$ are the intrinsic curvilinear coordinates of the membrane,
\begin{equation*}
dS = \sqrt{a} \,ds^1\,ds^2,  \qquad a = \textrm{det}\left(\frac{\partial x}{\partial s^\alpha} \cdot \frac{\partial x}{\partial s^\beta}\right)_{1 \leq \,\alpha, \,\beta\, \leq 2} \qquad \forall x \in \Omega.
\end{equation*}
We consider the following electrostatic potential energy of the protein-membrane system,
\begin{equation}\label{eq:G_gen}
G[\Gamma; \phi] = \DD\int_\Omega \left[-\frac{\varepsilon}{2}|\nabla\phi|^2 + f\phi - \chi_sB(\phi)\right]d\Omega + \frac{C}{\beta} \ln\left(\int_\Gamma \gamma(\phi) \; dS \Bigg/ \int_\Gamma \; dS\right).
\end{equation}
The electrostatic potential energy (also referred to as the electrostatic free energy) is determined by the positions of every interface $\Gamma$ and is expressed through the electrostatic potential $\phi: \Omega \to \reals$.  The electrostatic potential is designated by $\phi^s$, $\phi^p$, and $\phi^m$ in the domains $\Omega_s$, $\Omega_p$, and $\Omega_m$, respectively.  
The dielectric coefficient is $\varepsilon:\Omega \to \reals$. Distinct dielectric permittivities are defined 
in $\Omega_m$, $\Omega_s$, and $\Omega_p$, so $\varepsilon$ is a piecewise constant function,
\begin{equation}
\varepsilon(x) = \left\{\begin{split}
&\varepsilon_s && \qquad \textrm{if }x \in \Omega_s, \\
&\varepsilon_p && \qquad \textrm{if }x \in \Omega_p, \\
&\varepsilon_m && \qquad \textrm{if }x \in \Omega_m. \\
\end{split}\right.
\end{equation}
The fixed charge density inside the protein is given by $f:\Omega_p \to \reals$, an integrable function approximating the point charges in the protein.  The function $\chi_s$ denotes the characteristic function on $\Omega_s$ that is 1 on $\Omega_s$ and 0 elsewhere.  The function $B:\Omega \to \reals$ 
describes the electrostatic energy due to the mobile ions and is defined by
\begin{equation}\label{eq:B}
B(\phi) = \beta^{-1} \DD\sum_{j=1}^M c_j^\infty (e^{-\beta q_j \phi}-1), 
\end{equation}
where $\beta=1/(k_BT)$ is the inverse thermal energy, $M$ is the number of ionic species in the solvent, and $q_j$ and $c_j^\infty$ are the 
charge and bulk concentration of the $j^{th}$ ionic species, respectively \cite{Li2011}.  
The distribution of charged lipids $\rho[\Gamma]$ on the membrane is defined through the function $\gamma:\Omega \to \reals$ and satisfies the form
\begin{equation}\label{eq:rho_gen}
\rho[\Gamma] = \frac{C \gamma'(\phi)}{\beta q_l\DD\int_\Gamma \gamma(\phi) \; dS},
\end{equation}  
where $q_l$ is the charge of an individual lipid and $C$ is a dimensionless constant relate to the total quantity of charged lipids on membrane surfaces.  Here, we assume there is only one species of diffusive charged lipids in the membrane, but the model can easily be generalized to multiple lipid 
species. Note that \eqref{eq:rho_gen} is defined on the membrane boundaries $\Gamma_c$ and $\Gamma_e$ but not on the protein boundary $\Gamma_p$, as the lipids occupy $\Omega_m$ only.  

The main result of this article is the calculation of the dielectric boundary force induced by the protein on the charged membrane, which is computed through the shape derivative of the electrostatic potential energy, $\delta_\Gamma G[\Gamma]$.  The result is   
\begin{equation}\label{eq:boundary_force_copy}\begin{split}
F_n &= - \DD\frac{\varepsilon_s}{2}|\nabla\phi^s|^2 + \frac{\varepsilon_m}{2}|\nabla \phi^m|^2 - \varepsilon_m|\nabla\phi^m\cdot n|^2 + \varepsilon_m(\nabla \phi^s \cdot n)(\nabla \phi^m\cdot n)  \\
& \qquad  - B(\phi) - q_l\rho[\Gamma](\nabla\phi^s \cdot n).
\end{split}\end{equation}
where $\phi$ is the maximizer of the electrostatic potential energy $G[\Gamma; \cdot]$.  

The rest of the article is organized as follows.  
In Section \ref{sec:model}, we introduce the governing (Poisson-Boltzmann) equation for the protein-membrane electrostatic interactions.
In Section \ref{sec:shape_derivative}, we compute the shape derivative of the electrostatic energy functional, emphasizing the
treatment of the membrane surface charge distribution and the surface transformation.  Our main result, the dielectric boundary force on a charged surface, will follow from this computation.  
Finally, in Section \ref{sec:equivalence_to_MST}, we affirm that the force is physical by showing that it matches the divergence of the jump of Maxwell stress tensor. 
%
%
\section{Modeling electrostatics and charge distributions on membrane surfaces}\label{sec:model}

We now introduce the nonlinear Poisson-Boltzmann equation as
\begin{equation}\label{eq:PB}
\left\{\begin{gathered}\begin{aligned}
&\nabla\cdot(\varepsilon \nabla \phi) - \chi_s B'(\phi) = -f && \qquad \textrm{in }\Omega, \\
&[\phi] = 0 && \qquad\textrm{on }\Gamma_c,\Gamma_e,\Gamma_p, \\
&\varepsilon_s \frac{\partial \phi^s}{\partial n} = \varepsilon_m \frac{\partial \phi^m}{\partial n} - \rho[\Gamma] q_l && \qquad\textrm{on }\Gamma_c, \Gamma_e, \\
&\varepsilon_s \frac{\partial \phi^s}{\partial n} = \varepsilon_p \frac{\partial \phi^p}{\partial n} && \qquad\textrm{on }\Gamma_p,\\
&\phi = g && \qquad \textrm{on }\partial \Omega.
\end{aligned}\end{gathered}\right.
\end{equation}
This equation includes a general distribution of charged lipids which follows the form of \eqref{eq:rho_gen}.  The proceeding theorem establishes the relation between the electrostatic potential energy $G[\Gamma]$ for the lipid membrane boundaries $\Gamma_e$ and $\Gamma_c$ and the weak form of the nonlinear Poisson-Boltzmann equation.
In the theorem, we use the following notation for the Hilbert space,
\begin{equation}\label{def:H1g}
H^1_g(\Omega) = \{ \phi \in H^1(\Omega) | \phi = g \textrm{ on } \partial\Omega \}.
\end{equation}

\begin{theorem}\label{th:G_gen}
Let the electrostatic potential energy $G$ be
\begin{equation*}
G[\Gamma; \phi] = \DD\int_\Omega \left[-\frac{\varepsilon}{2}|\nabla\phi|^2 + f\phi - \chi_sB(\phi)\right]d\Omega + \frac{C}{\beta} \ln\left(\int_\Gamma \gamma(\phi) \; dS \Bigg/ \int_\Gamma \; dS\right),
\end{equation*}
where $\Gamma$ is any smooth boundary in $\Omega$.  Then there exists a critical point $\psi_0 \in H^1_g(\Omega)$ of $G[\Gamma; \cdot]$ which is also the weak solution of the nonlinear Poisson-Boltzmann equation \eqref{eq:PB}.
\end{theorem}
\begin{proof}
Let $\psi \in H^1_0(\Omega)$.  Consider the variation of the electrostatic potential energy:
\begin{equation*}\begin{split}
\left.\frac{dG[\Gamma;\phi+t\psi]}{dt}\right|_{t=0} \hspace{-4em}&\\
	&= \left.\frac{d}{dt}\left( \DD\int_\Omega \left[-\frac{\varepsilon}{2}|\nabla(\phi+t\psi)|^2  - \chi_sB(\phi+t\psi) + f(\phi+t\psi) \right]d\Omega \right. \right.  \\ 
	& \quad \qquad \left.\left. + \frac{C}{\beta} \ln\left(\int_\Gamma \gamma(\phi+t\psi) \; dS \Bigg/ \int_\Gamma \; dS\right) \right) \right|_{t=0}
\end{split}\end{equation*}
\begin{equation*}\begin{split}
\phantom{\left.\frac{dG[\Gamma;\phi+t\psi]}{dt}\right|_{t=0} \hspace{-4em}} &= \int_\Omega -\varepsilon \nabla(\phi+t\psi)\cdot (\nabla\psi) - \chi_s B'(\phi+t\psi)\psi + f(\phi+t\psi) \; d\Omega \\
	& \quad \qquad  \left. + \frac{C}{\beta} \left( \frac{\DD\int_\Gamma \gamma'(\phi+t\psi)\psi\; dS}{\DD\int_\Gamma \gamma(\phi+t\psi)\; dS}   \right) \right|_{t=0}  \\
\end{split}\end{equation*}
\begin{equation*}\begin{split}
\phantom{\left.\frac{dG[\Gamma;\phi+t\psi]}{dt}\right|_{t=0} \hspace{-4em}} 
	&= \DD\int_\Omega - \varepsilon (\nabla\phi)\cdot(\nabla\psi) - \chi_s B'(\phi)\psi + f\psi\; d\Omega + \int_\Gamma\frac{C\gamma'(\phi)\psi}{\beta \int_\Gamma \gamma(\phi) \; dS}\\
	&= \DD\int_\Omega - \varepsilon (\nabla\phi)\cdot(\nabla\psi) - \chi_s B'(\phi)\psi + f\psi \; d\Omega +\DD\int_\Gamma \rho[\Gamma] q_l \psi \; dS.
\end{split}\end{equation*}
%
%
Critical points of the electrostatic potential energy are obtained by setting the result above to 0,
\begin{equation}\label{eq:PB_var2}\begin{split}
\left(\int_{\Omega} -\varepsilon (\nabla \phi)\cdot(\nabla \psi) \; dX - 
\int_{\Omega_s} B'(\phi)\psi \;dX + \int_\Omega f\psi \;dX \right) &  \\
	+ \int_{\Gamma_c} \rho[\Gamma_c]q_l\psi\; dS +\int_{\Gamma_e} \rho[\Gamma_e]q_l\psi \;dS & = 0. 
\end{split}\end{equation}
Standard calculations show \eqref{eq:PB_var2} has a solution $\phi = \psi_0$ and that this critical point $\psi_0$ is the unique maximizer of $G[\Gamma; \cdot]$.  

Next, we will show that \eqref{eq:PB_var2} is equivalent to the weak form of \eqref{eq:PB}.
The weak form of the Poisson-Boltzmann equation for an arbitrary test function $\psi:\Omega\to\reals$ is 
\begin{equation}\label{eq:PB_variational}
\int_\Omega \nabla\cdot(\varepsilon \nabla \phi)\psi \; dX - \int_\Omega \chi_s B'(\phi)\psi \;dX = \int_\Omega -f\psi \;dX.
\end{equation}
Splitting the first integral in the domains $\Omega_s$, $\Omega_m$, and $\Omega_p$, and using the product rule gives 
%
%
%
\begin{equation*}\begin{split}
&\int_{\Omega_s}\nabla\cdot(\varepsilon_s \nabla \phi^s\psi) \; dX - \int_{\Omega_s}\varepsilon_s (\nabla \phi^s)\cdot(\nabla \psi) \; dX + \int_{\Omega_m}\nabla\cdot(\varepsilon_m \nabla \phi^m \psi) \; dX  \\
&\qquad - \int_{\Omega_m} \varepsilon_m (\nabla \phi^m)\cdot(\nabla\psi) \; dX + \int_{\Omega_p}\nabla\cdot(\varepsilon_p \nabla \phi^p \psi) \; dX - \int_{\Omega_p} \varepsilon_p (\nabla \phi^p)\cdot(\nabla\psi) \; dX \\
&\qquad - \int_{\Omega_s} B'(\phi)\psi \;dX = \int_\Omega -f\psi \;dX.
\end{split}\end{equation*}
Now combine the appropriate 
integrals 
and use the divergence theorem 
to get 
\begin{equation*}\begin{gathered}\begin{aligned}
&-\int_{\Omega} \varepsilon (\nabla \phi) \cdot (\nabla \psi) \; dX + \int_{\partial\Omega}\varepsilon_s \psi(\nabla \phi^s)\cdot n \; dS -\int_{\Gamma_c}\varepsilon_s \psi(\nabla \phi^s)\cdot n \; dS  \\
&\qquad - \int_{\Gamma_e}\varepsilon_s \psi(\nabla \phi^s)\cdot n \; dS -\int_{\Gamma_p}\varepsilon_s \psi(\nabla \phi^s)\cdot n \; dS + \int_{\Gamma_e} \varepsilon_m \psi(\nabla \phi^m)\cdot n \; dS  \\
&\qquad + \int_{\Gamma_c}\varepsilon_m \psi(\nabla \phi^m)\cdot n \; dS + \int_{\Gamma_p}\varepsilon_p \psi(\nabla \phi^p)\cdot n \; dS - \int_{\Omega_s} B'(\phi)\psi \;dX = -\int_\Omega f\psi \;dX.
\end{aligned}\end{gathered}\end{equation*}
Since 
$\psi$ is compactly supported on $\Omega$, the boundary integrals over $\partial \Omega$ are zero.  Combine the boundary integrals over $\Gamma$ and obtain
\begin{equation*}\begin{split}
&\left(-\int_{\Omega} \varepsilon (\nabla \phi)\cdot(\nabla \psi) \; dX  - \int_{\Omega_s} B'(\phi)\psi \;dX +\int_\Omega f\psi \;dX \right)  \\
& \qquad + \left(\int_{\Gamma_c} \psi(\varepsilon_m \nabla \phi^m - \varepsilon_s \nabla \phi^s)\cdot n \; dS \right) + \left(\int_{\Gamma_e} \psi(\varepsilon_m \nabla \phi^m - \varepsilon_s \nabla \phi^s)\cdot n \; dS \right) \\
& \qquad +  \left(\int_{\Gamma_p} \psi(\varepsilon_p \nabla \phi^p - \varepsilon_s \nabla \phi^s)\cdot n \; dS \right) = 0.
\end{split}\end{equation*}
Next, apply the boundary conditions 
to obtain \eqref{eq:PB_var2}, the variational form of the electrostatic potential energy.  Therefore, the maximizer $\psi_0$ of $G[\Gamma;\cdot]$ is also the weak solution to the nonlinear Poisson-Boltzmann equation. 
\end{proof}

Define the maximization of the electrostatic potential energy by
\begin{equation}\label{eq:G_max}
G[\Gamma] = \max_{\phi \in H^1_g(\Omega)} G[\Gamma;\phi] = G[\Gamma; \psi_0].
\end{equation}


We note that if $\gamma(\phi) = -e^{-q_l\beta\phi}$ and $C = \int_\Gamma \rho \; dS$ is the total number of charged lipids on $\Gamma$, a conserved quantity,
%
%
then the lipids follow the surface electrodiffusion equation at the steady state,
\begin{equation}\label{eq:drift-diffusion}
\frac{\partial \rho}{\partial t} = \nabla_s \cdot (D\nabla_s \rho + D q_l \rho \nabla_s \phi), 
\end{equation}
where $t$ is time, $D$ is the diffusion coefficient, and $\nabla_s, \nabla_s \cdot$ are the surface gradient and surface divergence operators, respectively.  
If we assume the lipid distribution $\rho[\Gamma]$ is governed by the surface electrodiffusion 
equation, then $\rho[\Gamma]$ follow the ideal Boltzmann distribution.
However, the practical distribution of lipids is subject to various constraints such as finite sizes and entropy 
conditions and may not follow this surface electrodiffusion equation \cite{KhelashviliG2009a,KiselevV2011a,ZhouY2012b}. 


\section{Shape derivative of electrostatic potential energy}  \label{sec:shape_derivative}

In this section, the main result $\delta_\Gamma G[\Gamma]$ is established through the method of shape derivatives.  For the shape derivative calculation, a velocity of the membrane movement is defined.  According to the Hadamard-Zol{\'e}sio structure theorem of shape calculus \cite{VanderZeePhDThesis}, the variation of the electrostatic potential energy \eqref{eq:G_gen} with respect to the position of the smooth boundary $\Gamma$, that is, the shape derivative, must be given by the inner product of the force with the normal component of the velocity of the deformation.  Thus, computing the shape derivative $\delta_\Gamma G[\Gamma]$ provides a formula for the force on a membrane $\Gamma$. 


\subsection{Velocity of the surface}

Define the velocity function $V\in C^\infty(\reals^3, \reals^3)$ by the following dynamical system,
\begin{equation}\label{eq:velocity}
\left\{\begin{gathered}\begin{aligned}
&\frac{dx}{dt} = V(x) \qquad \forall t >0,\\
&x(0) = X,
\end{aligned}\end{gathered}\right.
\end{equation}
where $X$ is the original position of the membrane and $x$ is the transformed position. Assume that $V$ is compactly supported near the bilayer membrane 
surfaces, i.e., $V(x) = 0$ if dist($x, \Gamma)>d$ for some $d>0$ where $\Gamma$ is either $\Gamma_c$ or $\Gamma_e$ and 
\begin{equation}\label{eq:dist}
d < \frac{1}{2} \min \;[\textrm{dist}\;(\Gamma_e, \partial\Omega), \textrm{dist}\;(\Gamma_c,\textrm{supp}\;(f)), \textrm{dist}\;(\Gamma_c,\textrm{supp}\;(f))].
\end{equation}
The condition in \eqref{eq:dist} prevents the exterior membrane surface from stretching beyond the containment domain $\Omega$ and also
prevents the either membrane surface from overlapping the center of an atom contributing to the charge density $f$ within the protein.  

The solution to \eqref{eq:velocity}  defines a diffeomorphism $T_t:\reals^3 \rightarrow \reals^3$ where $T_t(X) = x(t,X)$ maps the old coordinates $X$ into the transformed coordinates $x$.  By a Taylor expansion, we can approximate the map by 
\begin{align}
T_t(X) &= x(t,X)\notag\\
	&= x(0,X) + t\partial_t x(0,X) + \mathcal{O}(t^2)\notag\\
	&= X+tV(x(0,X)) + \mathcal{O}(t^2)\notag\\
	&= X +tV(X) + \mathcal{O}(t^2), \label{eq:transformation}
\end{align}
so that the map $T_t(X)$ agrees with the perturbation of the identity up to the leading term.  

The configuration under the transformation $T_t(X)$ results in a new electrostatic potential energy,  
\begin{equation} \label{eq:Gt}
G[\Gamma_t; \phi] = \DD\int_{\Omega_t} \left[-\frac{\varepsilon}{2}|\nabla\phi|^2 + f\phi - \chi_{s}B(\phi)\right]dx + \frac{C}{\beta} \ln\left(\int_{\Gamma_t} \gamma(\phi) \; dS \Bigg/ \int_{\Gamma_t} dS\right),
\end{equation}
where each of the functions are computed over the transformed regions $\Omega_t = T_t(\Omega)$, $\Gamma_t = T_t(\Gamma)$.  
Similar to the standard nonlinear Poisson-Boltzmann equation without surface charge distributions \cite{Li2011}, 
there is a unique maximizer $\psi_t \in H^1_g(\Omega) \cap L^\infty(\Omega)$ that maximizes \eqref{eq:Gt} over $H^1_g(\Omega)$. The maximum is
\begin{equation}\label{eq:Gt_max}
G[\Gamma_t] = \max_{\phi \in H^1_g(\Omega)} G[\Gamma_t;\phi] = G[\Gamma_t; \psi_t].
\end{equation}
Following the ideas of Theorem \ref{th:G_gen}, the same $\psi_t$ is also the unique weak solution to the transformed boundary value problem of the 
Poisson-Boltzmann equation,
\begin{equation}\label{eq:PB_t}
\left\{\begin{gathered}\begin{aligned}
&\nabla\cdot(\varepsilon \nabla \phi) - \chi_s B'(\phi) = -f && \qquad \textrm{in }\Omega_t, \\
&[\phi] = 0 && \qquad\textrm{on }\Gamma_{c_t},\Gamma_{e_t},\Gamma_{p_t}, \\
&\varepsilon_s \frac{\partial \phi^s}{\partial n} = \varepsilon_m \frac{\partial \phi^m}{\partial n} - 
\rho[\Gamma_t] q_l && \qquad\textrm{on }\Gamma_{c_t}, \Gamma_{e_t}, \\
&\varepsilon_s \frac{\partial \phi^s}{\partial n} = \varepsilon_p \frac{\partial \phi^p}{\partial n} && \qquad\textrm{on }\Gamma_{p_t}, \\
&\phi = g && \qquad \textrm{on }\partial \Omega.
\end{aligned}\end{gathered}\right.
\end{equation}


\subsection{Transformation properties}

The transformation $T_t(X)$ defined by \eqref{eq:velocity} acts on volumes and surfaces in $\Omega$.  Some useful properties of the transformation are outlined below. 

\subsubsection{Properties of the volume transformation}

The following are assumed properties of the transformation $T_t(X)$ defined by \eqref{eq:velocity} on a volume element $dX \in \reals^3$.  These properties will be used in the computation of the shape derivative and their justifications are found in \cite{Delfour2011}, among other sources.

\begin{enumerate}[(T1)]
\item  Let $X \in \mathbb{R}^3$ and $t\geq 0$.  Let $\nabla T_t(X)$ be the Jacobian matrix of $T_t$ at $X$ defined by $(\nabla T_t(X))_{ij} = \partial_jT_t^i(X)$, where $T_t^i$ is the $i$th component of $T_t$ ($i=1, 2, 3)$.  Let
\begin{equation}\label{eq:Jt}
J_t(X) = \textrm{det}\nabla T_t(X).
\end{equation}
For each $X$, the function $t\mapsto J_t(X)$ is in $C^{\infty}$ and at $X$,
\begin{equation}\label{eq:d_Jt}
\frac{dJ_t}{dt} = J_t(\nabla \cdot V) \circ T_t.
\end{equation}
At $t=0$, since no time has passed, $\nabla T_0 = I$ for any $x$ and so $J_0(X) = 1$.  The continuity of $J_t$ at $t=0$ implies $J_t >0 $ for $t>0$ small enough.
\item Define $A(t):\Omega \rightarrow \mathbb{R}$ for $t\geq 0$ small enough by
\begin{equation}\label{eq:A}
A(t) = J_t(\nabla T_t)^{-1}(\nabla T_t)^{-T}.
\end{equation}
where the notation $(\cdot)^{-T}$ denotes the inverse transpose.  
At each point in $\Omega$, 
\begin{equation}\begin{split}\label{eq:dA}
A'(t) &= \left[\left(\left(\nabla\cdot V\right)\circ T_t\right) - \left(\nabla T_t\right)^{-1}\left(\left(\nabla V\right)\circ T_t\right)\nabla T_t \right.\\
&\qquad - \left.\left(\nabla T_t\right)^{-1}\left(\left(\nabla V\right)\circ T_t\right)^{T}\left(\nabla T_t\right)\right]A(t).
\end{split}\end{equation}
%
\end{enumerate}


\subsubsection{Properties of the surface transformation}

In this section, we compute a useful property of the transformation $T_t$ on a surface element $dS \in \reals^2$.  Define the Jacobian of the surface transformation by
\begin{equation}\label{eq:J_s}
J_s(X,t) = (\det \nabla T_t(X))\left|\nabla T_t^{-T}  n(X)\right|
\end{equation}
as in \cite{Ciarlet1994}, where $n(X)$ is the normal unit normal vector at $X$.  Note that at time $t=0$, $J_s(X,0) = (\det I)\left|I  n\right| = 1\,$. 
  Analogous to the transformation on a volume element, the differential surface element is transformed by $ds = J_s \, dS$.  
The surface transformation property we wish to establish is stated as the following theorem: 
\begin{theorem}\label{th:S1}
The time derivative of the Jacobian of the surface transformation $J_s$ at $t=0$ is given by the surface divergence of the velocity, 
\begin{equation}\label{eq:dJ_s(0)}
\left.\frac{dJ_s}{dt}\right|_{t=0} = \nabla_s \cdot V .
\end{equation}
\end{theorem}
An elementary proof of this theorem is given in Appendix \ref{appa}. A similar result is given in \cite{Simon1980} (Lemma 4.2).  

Finally, a useful property on enclosed surfaces is shown.  
\begin{lemma}\label{lem:incompressible}
The surface divergence of any continuous function $F$ on an enclosed surface is 0, 
\begin{equation}\label{eq:incompressible}
\int_\Gamma (\nabla_s \cdot F) \; dS = 0.
\end{equation}
\end{lemma}
\begin{proof}
%
%
%
The proof follows directly from the divergence theorem and the fact that the measure of the boundary $\partial \Gamma$ is zero.  
\end{proof}


\subsection{Shape derivative calculation}

The central theorem of the paper establishes the shape derivative of the electrostatic potential energy to the boundary, 
which provides a convenient way to extract $F_n$, the normal component of the dielectric boundary force,  
\begin{equation}\label{eq:F_def}
\delta_{\Gamma}G[\Gamma] =  \int_\Gamma -F_n (V \cdot n) \; dS .
\end{equation}
This is the Hadamard-Zol\'{e}sio structure theorem \cite{VanderZeePhDThesis, Stratton1941}.  We point out that only the normal component of the force determines the motion on the boundary.  The following theorem shows that this force is 
\begin{equation}\label{eq:boundary_force}\begin{split}
F_n &= -\DD\frac{\varepsilon_s}{2}|\nabla\psi_0^s|^2 + \frac{\varepsilon_m}{2}|\nabla \psi_0^m|^2 - \varepsilon_m|\nabla\psi_0^m\cdot n|^2 + \varepsilon_m(\nabla \psi_0^s \cdot n)(\nabla \psi_0^m\cdot n)  \\
& \qquad  - B(\psi_0) - q_l\rho[\Gamma](\nabla\psi_0^s \cdot n),
\end{split}\end{equation}
which is our main result.

Let $V \in C^\infty(\mathbb{R}^3,\mathbb{R}^3)$ be a smooth map that vanishes outside a small neighborhood of the membrane surface $\Gamma$.  That is, $V(X) = 0$ if dist$(X,\Gamma) >d$ for some $d>0$ satisfying \eqref{eq:dist}.  Let the transformations $T_t (t\geq 0)$ be defined by \eqref{eq:velocity}.  
For $t>0$ the electrostatic free energy is given by \eqref{eq:Gt_max}, where the functional $G[\Gamma_t;\cdot]$ is given in \eqref{eq:Gt} and $\psi_t$ is the weak solution to \eqref{eq:PB_t}.  
For $t=0$, the electrostatic free energy is given by \eqref{eq:G_max}, where the functional $G[\Gamma; \cdot]$ is given in \eqref{eq:G_gen} and $\psi_0$ is the weak solution to \eqref{eq:PB}.

\begin{theorem}\label{th:force} Assume $f \in H^1(\Omega)$.  Then the shape derivative of the electrostatic free 
energy $G[\Gamma]$ in the direction of $V$ is given by 
\begin{align}
&\delta_{\Gamma} G[\Gamma] = \DD\int_{\Gamma} \Bigg(\frac{\varepsilon_s}{2}|\nabla\psi_0^s|^2 - \frac{\varepsilon_m}{2}|\nabla \psi_0^m|^2 + \varepsilon_m|\nabla\psi_0^m\cdot n|^2 - \varepsilon_m(\nabla \psi_0^s \cdot n)(\nabla \psi_0^m\cdot n)  \nonumber \\
& \qquad \qquad \qquad  + B(\psi_0) + q_l\rho[\Gamma](\nabla\psi_0^s \cdot n)  \Bigg) (V \cdot n) \; dS. \label{eq:shape_deriv}
\end{align}
This shape derivative is computed as the variation of the electrostatic potential energy with respect to the membrane surface. 
\end{theorem}
\begin{proof} The proof is divided into four steps.  

\begin{enumerate}[i.]
\item First, the energy functional is computed in the transformed coordinates through a new function $z(t,\phi)$.  
A change of variables brings $z$ back to the reference coordinates, and then $z$ is differentiated with respect to time.  
\item Second, the difference quotient corresponding to the shape derivative is squeezed between two realizations of $\partial_t z$.  
\item Third, the inequality is passed to the limit as $t \to 0$ and we will show that the two realizations of $\partial_t z$ are identical in the limit, and hence equal to the shape derivative.  
\item Fourth and finally, the result is simplified to match the final form. 
\end{enumerate}

The computations to determine the shape derivative of $\eqref{eq:G_gen}$ is completed in two calculations by splitting $\eqref{eq:G_gen}$ into two components,
\begin{equation}\label{eq:G1}
G_1[\Gamma; \phi] = \DD\int_\Omega \left[-\frac{\varepsilon}{2}|\nabla\phi|^2 + f\phi - \chi_sB(\phi)\right]d\Omega 
\end{equation}
and
\begin{equation}\label{eq:G2}
G_2[\Gamma; \phi] = \frac{C}{\beta} \ln\left(\int_\Gamma \gamma(\phi) \; dS \Bigg/ \int_\Gamma \; dS\right),
\end{equation}
where $G[\Gamma;\phi] = G_1[\Gamma;\phi] + G_2[\Gamma;\phi]$. The details of the computations for $G_1[\Gamma;\phi]$ are 
omitted as they appear in \cite{Li2011}.  

\noindent \textbf{Step 1}.  Let $t\geq 0$ be sufficiently small.

Note that $\phi \in H^1_g(\Omega) \to \phi \circ T_t^{-1} \in H^1_g(\Omega)$ is an isomorphism, so by \eqref{eq:Gt_max}, we have
$$
G[\Gamma_t] = \max_{\phi \in H^1_g(\Omega)} G[\Gamma_t; \phi \circ T_t^{-1}].
$$
This maps the transformed coordinates $x= T_t(X)$ back into the original coordinates, for which information is known.  Let $\phi \in H^1(\Omega) \cap L^\infty(\Omega)$ and $t\geq0$ and denote 
\begin{equation}\label{eq:z}
z(t,\phi) = G[\Gamma_t; \phi \circ T_t^{-1}].
\end{equation}
The function $z$ is also split into two components, $z_1$, corresponding to $G_1$, and $z_2$ corresponding to $G_2$:
\begin{align}
&z_1(t,\phi) = G_1[\Gamma_t; \phi \circ T_t^{-1}], \label{eq:z1}\\
&z_2(t,\phi) = G_2[\Gamma_t; \phi \circ T_t^{-1}], \label{eq:z2}
\end{align}
so that $z(t,\phi) = z_1(t,\phi) + z_2(t,\phi)$.  
Since the details of $z_1(t,\phi)$ appear exactly in \cite{Li2011}, we only consider $z_2(t,\phi)$ here.
Using \eqref{eq:z2}, \eqref{eq:G2}, and a transformation of coordinates $x=T_t(X)$,
\begin{align}
z_2(t, \phi) &= G_2[\Gamma_t;\phi \circ T_t^{-1}] \notag\\
	&= \frac{C}{\beta}\ln\left(\int_{\Gamma_t} \gamma(\phi \circ T_t^{-1})(x) \; dS(x) \Bigg/  \int_{\Gamma_t} \; dS(x)\right)  \notag\\
	&= \frac{C}{\beta}\ln\left(\DD\int_{\Gamma_0} \gamma( \phi (X)) J_s(X,t) \; dS(X) \Bigg/  \DD\int_{\Gamma_0} J_s(X,t) \; 
dS(X)\right), \label{eq:z2_2}
\end{align}
where the formula for $J_s$ is given by \eqref{eq:J_s}.  Note that $\phi(X)$ 
does not depend on $t$.  Differentiating with respect to $t$,
\begin{equation}\label{eq:dz2dt}\begin{gathered}\begin{aligned}
\partial_t z_2(t,\phi) &= \frac{C}{\beta} \left( \frac{\DD\int_{\Gamma_0} \gamma(\phi) \frac{dJ_s(X,t)}{dt} \; dS(X)}{\DD\int_{\Gamma_0} \gamma(\phi) J_s(X,t) \; dS(X) } - \frac{\DD\int_{\Gamma_0} \frac{dJ_s(X,t)}{dt} \; dS(X)}{\DD\int_{\Gamma_0} J_s(X,t) \; dS(X) } \right).
\end{aligned}\end{gathered}\end{equation}
Then, the full form of $\partial_t z(t,\phi)$ is given by combining the calculation in \cite{Li2011} and \eqref{eq:dz2dt},
\begin{align}
\partial_t z(t,\phi) =  & \DD\int_\Omega \left[-\frac{\varepsilon}{2}A'(t)\nabla\phi\cdot\nabla\phi + ((\nabla \cdot (fV))\circ T_t)\phi J_t - \chi_sB(\phi)((\nabla \cdot V)\circ T_t) J_t \right]dX  \nonumber \\
	& + \frac{C}{\beta} \left( \frac{\DD\int_{\Gamma_0} \gamma( \phi) \frac{dJ_s(X,t)}{dt} \; dS(X)}{\DD\int_{\Gamma_0} \gamma(\phi) J_s(X,t) \; dS(X) } - \frac{\DD\int_{\Gamma_0} \frac{dJ_s(X,t)}{dt} \; dS(X)}{\DD\int_{\Gamma_0} J_s(X,t) \; dS(X) } \right), \label{eq:dzdt}
\end{align}
where $A'(t)$ is given by \eqref{eq:dA}.


\vspace{1em}

\textbf{Step 2.}  Let $t \in (0,\tau]$.  Since $\psi_t \in H^1_g(\Omega) \cap L^\infty (\Omega)$ and $\psi_0 \in H^1_g(\Omega) \cap L^\infty(\Omega)$ maximize $G[\Gamma_t; \cdot]$ and $G[\Gamma;\cdot]$ over $H_g^1(\Omega)$, respectively, we have
\begin{align}
G[\Gamma_t; \psi_0 \circ T_t^{-1}] \leq G[\Gamma_t; \psi_t] = G[\Gamma_t], \label{eq:squeeze1}\\
G[\Gamma; \psi_t] \leq G[\Gamma; \psi_0] = G[\Gamma], \label{eq:squeeze2}\\
G[\Gamma; \psi_t \circ T_t] \leq G[\Gamma; \psi_0] = G[\Gamma].\label{eq:squeeze3}
\end{align}

By \eqref{eq:squeeze1}, 
\begin{equation*}
\frac{G[\Gamma_t;\psi_0\circ T_t^{-1}]-G[\Gamma;\psi_0]}{t} \leq \frac{G[\Gamma_t] - G[\Gamma]}{t}.
\end{equation*}

By \eqref{eq:squeeze3},
\begin{equation*}
\frac{G[\Gamma_t] - G[\Gamma]}{t} \leq \frac{G[\Gamma_t; \psi_t] - G[\Gamma; \psi_t \circ T_t]}{t}.
\end{equation*}

Putting the two inequalities together with the definition in \eqref{eq:z} gives
\begin{equation}
\frac{z(t,\psi_0)-z(0,\psi_0)}{t} \leq \frac{G[\Gamma_t]-G[\Gamma]}{t} \leq \frac{z(t, \psi_t \circ T_t) - z(0, \psi_t \circ T_t)}{t}.
\end{equation}

Notice that the far left expression of the inequality is the secant line of $z(\cdot,\psi_0)$ from 0 to $t \leq \tau$, and the far right expression is the secant line of $z(\cdot,\psi_t\circ T_t)$ from 0 to $t\leq \tau$.  Since $z$ is differentiable on $t\geq 0$, we can apply the Mean Value Theorem to each expression.  That is, there exists $\xi(t), \eta(t) \in [0,t]$ for each $t\in (0,\tau]$ such that 

\begin{equation}\label{eq:sandwich}
\partial_t z(\xi(t),\psi_0) \leq \frac{G[\Gamma_t]-G[\Gamma]}{t} \leq \partial_t z (\eta(t),\psi_t\circ T_t) 
\qquad \forall t\in (0,\tau].
\end{equation}


\vspace{1em}

\textbf{Step 3.} We claim 
\begin{align}
&\lim_{t\to 0} \partial_t z(\xi(t), \psi_0) = \partial_t z(0, \psi_0), \label{eq:step3a}\\
&\lim_{t\to 0} \partial_t z(\eta(t), \psi_t\circ T_t) = \partial_t z(0,\psi_0). \label{eq:step3b}
\end{align}

As in Step 1, the proofs for \eqref{eq:step3a} and \eqref{eq:step3b} are split into the corresponding $z_1$ and $z_2$ components,
\begin{align}
&\lim_{t\to 0} \partial_t z_1(\xi(t), \psi_0) = \partial_t z_1(0, \psi_0),  \label{eq:step3a1}\\
&\lim_{t\to 0} \partial_t z_1(\eta(t), \psi_t\circ T_t) = \partial_t z_1(0,\psi_0), \label{eq:step3b1} \\
&\lim_{t\to 0} \partial_t z_2(\xi(t), \psi_0) = \partial_t z_2(0, \psi_0),  \label{eq:step3a2}\\
&\lim_{t\to 0} \partial_t z_2(\eta(t), \psi_t\circ T_t) = \partial_t z_2(0,\psi_0). \label{eq:step3b2}
\end{align}

The justifications of \eqref{eq:step3a1} and \eqref{eq:step3b1} are found in \cite{Li2011}, so only \eqref{eq:step3b2} is proven here, noting that \eqref{eq:step3a2} is similar.  Let $\eta(t) \in [0,t]$ and do consider passing $\partial_tz_2(\eta(t), \psi_t\circ T_t)$ to the limit as $t \to 0$,
\begin{equation*}\begin{gathered}\begin{aligned}
\partial_t z_2(\eta(t)&, \psi_t \circ T_t) = \resizebox{.75\hsize}{!}{$\DD\frac{C}{\beta} \left( \frac{\DD\int_{\Gamma_0} \gamma (\psi_t \circ T_t) \frac{dJ_s(X,\eta(t))}{dt} \; dS(X)}{\DD\int_{\Gamma_0} \gamma(\psi_t \circ T_t) J_s(X,\eta(t)) \; dS(X) } - \frac{\DD\int_{\Gamma_0} \frac{dJ_s(X,\eta(t))}{dt} \; dS(X)}{\DD\int_{\Gamma_0} J_s(X,\eta(t)) \; dS(X) } \right)$} \\
	&\to \resizebox{.65\hsize}{!}{$\DD\frac{C}{\beta} \left( \frac{\DD\int_{\Gamma_0} \gamma(\psi_0) \frac{dJ_s(X,0)}{dt} \; dS(X)}{\DD\int_{\Gamma_0} \gamma(\psi_0) J_s(X,0) \; dS(X) } - \frac{\DD\int_{\Gamma_0} \frac{dJ_s(X,0)}{dt} \; dS(X)}{\DD\int_{\Gamma_0} J_s(X,0) \; dS(X) } \right)$}. 
\end{aligned}\end{gathered}\end{equation*}
By passing \eqref{eq:sandwich} to the limit as $t\to0$ and using \eqref{eq:step3a} and \eqref{eq:step3b} through the computation \eqref{eq:dzdt} shows that
\begin{align}
\delta_{\Gamma} G[\Gamma] &= \lim_{t\to0} \frac{G[\Gamma_t] - G[\Gamma]}{t} = \partial_tz(0,\psi_0) \notag \\
	&= \DD\int_\Omega \left[-\frac{\varepsilon}{2}A'(0)\nabla\psi_0\cdot\nabla\psi_0 + (\nabla \cdot (fV))\psi_0 - \chi_s B(\psi_0)(\nabla \cdot V) \right]dX \notag\\
	& \qquad + \frac{C}{\beta} \left( \frac{\DD\int_{\Gamma_0} \gamma(\psi_0) \frac{dJ_s(X,0)}{dt} \; dS(X)}{\DD\int_{\Gamma_0} \gamma( \psi_0) J_s(X,0) \; dS(X) } - \frac{\DD\int_{\Gamma_0} \frac{dJ_s(X,0)}{dt} \; dS(X)}{\DD\int_{\Gamma_0} J_s(X,0) \; dS(X) } \right). \label{eq:shape1}
\end{align}


\vspace{1em}
\textbf{Step 4.} 
The final step is to simplify the form of \eqref{eq:shape1} so that it matches the right-hand side of \eqref{eq:shape_deriv}.  
We begin our simplifications with $z_1$ and then move to $z_2$.  The combined results formulate the shape derivative. 
In \cite{Li2011}, $\partial_t z_1(0,\psi_0)$ is simplified to
\begin{equation}\begin{split}\label{eq:dz1dt_2}
\partial_t z_1(0,\psi_0)  = &
\DD\int_{\Gamma} \frac{\varepsilon_s}{2} |\nabla \psi_0^s|^2(V\cdot n) \; dS - \DD\int_{\Gamma} \frac{\varepsilon_m}{2} |\nabla \psi_0^m|^2(V\cdot n) \; dS  \\
& - \DD\int_{\Gamma} \varepsilon_s (\nabla \psi_0^s\cdot n)(V\cdot \nabla\psi_0^s) \; dS + \DD\int_{\Gamma} \varepsilon_m (\nabla \psi_0^m\cdot n)(V\cdot \nabla \psi_0^m) \; dS  \\
& + \DD\int_\Gamma B(\psi_0)(V\cdot n) \; dS. 
\end{split}\end{equation}
By the continuity of $\psi_0$ on $\Gamma$,
\begin{equation}\label{eq:above1}
\nabla(\psi_0^s -\psi_0^m) = (\nabla\psi_0^s\cdot n - \nabla \psi_0^m\cdot n)n \qquad\textrm{on } \Gamma.
\end{equation}
Now apply the new interface condition:
\begin{equation}\label{eq:new_jump}
\varepsilon_s\nabla \psi_0^s \cdot n = \varepsilon_m\nabla \psi_0^m\cdot n - q_l\rho[\Gamma],
\end{equation}
where $\rho[\Gamma]$ follows \eqref{eq:rho_gen}.  There are two equivalent approaches from here.  Consider the third and fourth terms of \eqref{eq:dz1dt_2}.  Substituting \eqref{eq:new_jump} into the third term and using \eqref{eq:above1} gives
\begin{align*}
& - \DD\int_{\Gamma} \varepsilon_s (\nabla \psi_0^s\cdot n)(V\cdot \nabla\psi_0^s) \; dS + \DD\int_{\Gamma} \varepsilon_m (\nabla \psi_0^m\cdot n)(V\cdot \nabla \psi_0^m) dS \\
& \quad = - \DD\int_{\Gamma} (\varepsilon_m\nabla \psi_0^m\cdot n - q_l\rho[\Gamma])(V\cdot \nabla\psi_0^s) \; dS + \DD\int_{\Gamma} \varepsilon_m (\nabla \psi_0^m\cdot n)(V\cdot \nabla \psi_0^m) \; dS \\
& \quad= \DD\int_{\Gamma} -(\varepsilon_m\nabla \psi_0^m\cdot n)V\cdot \nabla(\psi_0^s -\psi_0^m) + q_l\rho[\Gamma](V\cdot \nabla\psi_0^s) \; dS \\
& \quad= \DD\int_{\Gamma} -(\varepsilon_m\nabla \psi_0^m\cdot n)V\cdot (\nabla\psi_0^s \cdot n -\nabla\psi_0^m\cdot n)n + q_l\rho[\Gamma](V\cdot \nabla\psi_0^s) \; dS \\
& \quad= \DD\int_{\Gamma} \varepsilon_m \left|\nabla \psi_0^m\cdot n\right|^2 (V\cdot n) - \varepsilon_m(\nabla\psi_0^s\cdot n)(\nabla\psi_0^m\cdot n)(V\cdot n) + q_l\rho[\Gamma](V\cdot \nabla\psi_0^s) \; dS.
\end{align*}
Substituting this term back into \eqref{eq:dz1dt_2} and combining terms leads to the final result for $\partial_tz_1$,
\begin{align}
&\partial_t z_1(0,\psi_0) = \DD\int_{\Gamma} \frac{\varepsilon_s}{2} |\nabla \psi_0^s|^2(V\cdot n) \; dS - \DD\int_{\Gamma} \frac{\varepsilon_m}{2} |\nabla \psi_0^m|^2(V\cdot n) \; dS  \notag\\
&\qquad+ \resizebox{.9\hsize}{!}{$\DD\int_{\Gamma} \varepsilon_m \left|\nabla \psi_0^m\cdot n\right|^2 (V\cdot n) - \varepsilon_m(\nabla\psi_0^s\cdot n)(\nabla\psi_0^m\cdot n)(V\cdot n) + q_l\rho[\Gamma](V\cdot \nabla\psi_0^s) \; dS$} \notag\\
& \qquad + \DD\int_\Gamma B(\psi_0)(V\cdot n) \; dS \notag\\
&\quad = \DD\int_\Gamma \left[ \left(\frac{\varepsilon_s}{2}|\nabla\psi_0^s|^2 - \frac{\varepsilon_m}{2}|\nabla \psi_0^m|^2 + \varepsilon_m|\nabla\psi_0^m\cdot n|^2 - \varepsilon_m(\nabla \psi_0^s \cdot n)(\nabla \psi_0^m\cdot n) \right.\right. \notag\\
&\quad \quad \left.\left.\phantom{\frac{}{}}+B(\psi_0)\right)(V \cdot n) + q_l\rho[\Gamma](V\cdot \nabla \psi_0^s) \right]\; dS \label{eq:dzdt_final_b}
\end{align}
An alternative form of \eqref{eq:dzdt_final_b} can be obtained by substituting \eqref{eq:new_jump} and keeping $(\varepsilon_s\nabla\psi_0^s\cdot n)$ instead. 


Next, the shape derivative from $\partial_tz_2$ is simplified.  Evaluating \eqref{eq:dz2dt} at $t=0$ and $\phi = \psi_0$ gives,
\begin{align*}
\partial_tz_2(0,\psi_0) &= \frac{C}{\beta} \left( \frac{\DD\int_{\Gamma_0} \gamma(\psi_0) \frac{dJ_s(X,0)}{dt} \; dS(X)}{\DD\int_{\Gamma_0} \gamma(\psi_0) J_s(X,0) \; dS(X) } - \frac{\DD\int_{\Gamma_0} \frac{dJ_s(X,0)}{dt} \; dS(X)}{\DD\int_{\Gamma_0} J_s(X,0) \; dS(X) } \right). 
\end{align*}
The time derivative of $J_s$ at $t=0$ is the surface divergence of the velocity, 
as calculated in \eqref{eq:dJ_s(0)}. Recall that $J_s(X,0) = 1$. Using these two facts, the above equation is equivalent to   
\begin{align}\label{eq:shape_G2}
\partial_tz_2(0,\psi_0) &= \frac{C}{\beta} \left( \frac{\DD\int_{\Gamma_0} \gamma( \psi_0) (\nabla_s \cdot V ) \; dS(X)}{\DD\int_{\Gamma_0} 
\gamma( \psi_0)  \; dS(X) } - \frac{\DD\int_{\Gamma_0} (\nabla_s \cdot V ) \; dS(X)}{\DD\int_{\Gamma_0} \; dS(X) } \right).
\end{align}
By Lemma \ref{lem:incompressible}, since the velocity $V$ is a continuous function, the second term of \eqref{eq:shape_G2} is zero.
%
%
Applying the product rule to the 
remaining term,
\begin{align}\label{eq:shape_G2_3}
\partial_tz_2(0,\psi_0) &= \frac{C}{\beta} \left( \frac{\DD\int_{\Gamma_0} \nabla_s \cdot (V \gamma(\psi_0)) - V \cdot (\nabla_s \gamma(\psi_0)) \; dS(X)}{\DD\int_{\Gamma_0} \gamma(\psi_0)  \; dS(X) }  \right).
\end{align}
By applying Lemma \ref{lem:incompressible} again to \eqref{eq:shape_G2_3} to remove the surface divergence of the continuous function $V \gamma(\psi_0)$, and by expanding the surface gradient in the remaining term, we have
\begin{align}
\partial_tz_2(0,\psi_0) &= \frac{C}{\beta} \left( \frac{\DD\int_{\Gamma_0} - V \cdot (\gamma'(\psi_0)\nabla_s \psi_0 ) \; dS(X)}
{\DD\int_{\Gamma_0} \gamma(\psi_0)  \; dS(X) }  \right) \notag 
\\
%
%
&= -\DD\int_{\Gamma_0}  q_l \rho[\Gamma](V \cdot \nabla_s\psi_0)  \; dS(X).  \label{eq:shape_G2_final}
\end{align}
where the last step uses the definition of $\rho[\Gamma]$ in \eqref{eq:rho_gen} at $\phi = \psi_0$.  
%
%
Finally, combine \eqref{eq:dzdt_final_b} and \eqref{eq:shape_G2_final} for the full form of the shape derivative,
\begin{align*}
\delta_{\Gamma} G[\Gamma] &= \partial_t z_1 + \partial_t z_2 \\
& = \DD\int_{\Gamma_0} \left[ \left(\frac{\varepsilon_s}{2}|\nabla\psi_0^s|^2 - \frac{\varepsilon_m}{2}|\nabla \psi_0^m|^2 + \varepsilon_m|\nabla\psi_0^m\cdot n|^2 - \varepsilon_m(\nabla \psi_0^s \cdot n)(\nabla \psi_0^m\cdot n) \right.\right. \\
& \quad \left.\left.\phantom{\frac{}{}}+B(\psi_0)\right)(V \cdot n) + q_l\rho[\Gamma](V\cdot \nabla \psi_0^s) \right]\; dS - \int_{\Gamma_0}  q_l \rho[\Gamma] (V \cdot \nabla_s\psi_0)  \; dS.
\end{align*}
By continuity of $\psi_0$, $\psi_0^s = \psi_0^m$ on $\Gamma$
, and hence the last two terms reduce:
\begin{align*}
&\int_{\Gamma_0}  q_l\rho[\Gamma](V\cdot \nabla \psi_0^s) \; dS - \int_{\Gamma_0}  q_l \rho[\Gamma] (V \cdot \nabla_s\psi_0)  \; dS   \\
& \quad = \int_{\Gamma_0} q_l \rho[\Gamma] V \cdot (\nabla \psi_0^s - \nabla_s \psi_0^s) \; dS\\
& \quad = \int_{\Gamma_0} q_l \rho[\Gamma] V \cdot (\nabla \psi_0^s -[\nabla \psi_0^s - (\nabla \psi_0^s\cdot n)n])  \; dS\\
& \quad = \int_{\Gamma_0} q_l \rho[\Gamma] (\nabla\psi_0^s\cdot n)(V \cdot n) \; dS.
\end{align*}
%
%
%
Using this result, the shape derivative is simplified to the final result,
\begin{align}
&\delta_{\Gamma} G[\Gamma] = \DD\int_{\Gamma_0} \Bigg(\frac{\varepsilon_s}{2}|\nabla\psi_0^s|^2 - \frac{\varepsilon_m}{2}|\nabla \psi_0^m|^2 + \varepsilon_m|\nabla\psi_0^m\cdot n|^2 - \varepsilon_m(\nabla \psi_0^s \cdot n)(\nabla \psi_0^m\cdot n)  \nonumber \\
& \qquad \qquad \qquad  + B(\psi_0) + q_l\rho[\Gamma](\nabla\psi_0^s \cdot n)  \Bigg) (V \cdot n) \; dS. \label{eq:shape_final}
\end{align}

The alternative form of $\delta_{\Gamma}G[\Gamma]$ based on the alternative form of \eqref{eq:dzdt_final_b} is given by
\begin{align}
&\delta_{\Gamma}^{\rm(alt)} G[\Gamma] = \DD\int_{\Gamma_0} \Bigg(\frac{\varepsilon_s}{2}|\nabla\psi_0^s|^2 - \frac{\varepsilon_m}{2}|\nabla \psi_0^m|^2 - \varepsilon_s|\nabla\psi_0^s\cdot n|^2 + \varepsilon_s(\nabla \psi_0^s \cdot n)(\nabla \psi_0^m\cdot n)  \nonumber \\
& \qquad \qquad \qquad  + B(\psi_0) + q_l \rho[\Gamma](\nabla\psi_0^m \cdot n)\Bigg) (V \cdot n) \; dS. \label{eq:shape_final2}
\end{align}
Notice that the two forms are equivalent, since subtracting \eqref{eq:shape_final2} from \eqref{eq:shape_final} using the second interface condition in \eqref{eq:PB} results in zero.
%
%
\end{proof}


\section{Equivalence to the Maxwell stress tensor} \label{sec:equivalence_to_MST}
In classical electromechanics, the force on charged interfaces is given by the jump in the Maxwell stress tensor \cite{Gilson93, Landau1995, Sharp1990a, Stratton1941}.  The Maxwell stress tensor is defined by 
\begin{equation}\label{eq:MST}
M = \varepsilon E \otimes E - \frac{\varepsilon}{2}|E|^2I - \chi_s B(\psi)I, 
\end{equation}
where $E =-\nabla \psi$, not to be confused with the bending energy \eqref{eq:bending_energy}.
Since $\varepsilon$ and $\nabla \psi$ are discontinuous across the membrane interfaces $\Gamma_e$ and $\Gamma_c$, the Maxwell stress tensor $M$ is also discontinuous.  In the domains $\Omega_m$ and $\Omega_s$, we refer to the Maxwell stress tensor as $M^m$ and $M^s$, respectively.  
We verify that the force $F_n$ given by \eqref{eq:shape_deriv} is equivalent to the jump of the normal stress given by the
Maxwell stress tensor \cite{Cai2011}:
\begin{align}
F_n &= n\cdot M^sn-n\cdot M^mn \notag\\
	&= \Big(\varepsilon_s|\nabla\psi^s\cdot n|^2-\frac{\varepsilon_s}{2}|\nabla\psi^s|^2-\chi_sB(\psi^s)\Big) - \Big(\varepsilon_m|\nabla\psi^m\cdot n|^2-\frac{\varepsilon_m}{2}|\nabla\psi^m|^2 \Big) \notag\\
	&= -\frac{\varepsilon_s}{2}|\nabla\psi^s|^2 + \frac{\varepsilon_m}{2}|\nabla\psi^m|^2 + \varepsilon_s|\nabla\psi^s\cdot n|^2 - \varepsilon_m|\nabla\psi^m\cdot n|^2 -B(\psi)  \notag\\
	& = -\frac{\varepsilon_s}{2}|\nabla\psi^s|^2 + \frac{\varepsilon_m}{2}|\nabla\psi^m|^2 + \varepsilon_s(\nabla\psi^s\cdot n)(\nabla\psi^s\cdot n)- \varepsilon_m(\nabla\psi^m\cdot n)(\nabla\psi^s\cdot n) \notag\\
	& \qquad   +\varepsilon_m(\nabla\psi^s\cdot n)(\nabla\psi^m\cdot n)-\varepsilon_m(\nabla\psi^m\cdot n)(\nabla\psi^m\cdot n)  -B(\psi)\notag \\
	&= -\frac{\varepsilon_s}{2}|\nabla\psi^s|^2 + \frac{\varepsilon_m}{2}|\nabla\psi^m|^2 - \varepsilon_m|\nabla\psi^m\cdot n|^2 +\varepsilon_m(\nabla\psi^s\cdot n)(\nabla\psi^m\cdot n)  \notag \\
	& \qquad -B(\psi) +\Big(\varepsilon_s(\nabla\psi^s\cdot n)-\varepsilon_m(\nabla\psi^m\cdot n)\Big)(\nabla\psi^s\cdot n)  \notag \\
	&= -\frac{\varepsilon_s}{2}|\nabla\psi^s|^2 + \frac{\varepsilon_m}{2}|\nabla\psi^m|^2 - \varepsilon_m|\nabla\psi^m\cdot n|^2 +\varepsilon_m(\nabla\psi^s\cdot n)(\nabla\psi^m\cdot n)  \notag\\
	&\qquad  -B(\psi) - q_l\rho[\Gamma] (\nabla\psi^s\cdot n).
\end{align}
This matches \eqref{eq:boundary_force} exactly, suggesting that our definition of the electrostatic forces on
charged dielectric interfaces is physical and can reproduce the correct dynamics when coupled to the mechanical
models of bilayer membranes.


\section*{Acknowledgments}
The authors thank Benzhuo Lu and Bo Li for many helpful discussions.

\appendix


\section*{Appendix}\label{appa}

Proof of Theorem \ref{th:S1}.
\begin{proof}
Recall the formula for the Jacobian of the surface transformation, 
\begin{equation}\label{eq:J_s_app}
J_s(X,t) = (\det \nabla T_t(X))|\nabla T_t^{-T}  n(X)|.
\end{equation}
The derivative of \eqref{eq:J_s_app} is given by
\begin{equation}\label{eq:dJ_s}
\frac{dJ_s}{dt} = (J_t(\nabla \cdot V)\circ T_t)\left|\nabla T_t^{-T} n\right| 
	 + 
\frac{\det\nabla T_t}{\left|\nabla T_t^{-T}  n\right|}(\nabla T_t^{-T}  n)\cdot\frac{d}{dt}(\nabla T_t^{-T}  n).
\end{equation}
The above computation uses the product rule and \eqref{eq:d_Jt}.  We proceed by simplifying \eqref{eq:dJ_s}.  The derivative with respect to $t$ of $(\nabla T_t^{-T} n)$ is 
\begin{equation}\label{eq:dTndt}
\frac{d}{dt}(\nabla T_t^{-T}  n) = \frac{d(\nabla T_t^{-T})}{dt} n + \nabla T_t^{-T} \frac{dn}{dt}. 
\end{equation}
There is an alternative expression to \eqref{eq:dTndt} that will prove easier for future calculations. 
Let $\vec{a} = \nabla T_t^{-T} n$ so that $n = \nabla T_t^T \vec{a}$.  Then, computing the derivative of $n$,
\begin{equation*} \label{eq:dTn_dt}
\frac{dn}{dt} = \frac{d(\nabla T_t^T \vec{a})}{dt}
	= \frac{d(\nabla T_t^T)}{dt} \vec{a} + \nabla T_t^T \frac{d\vec{a}}{dt}.
\end{equation*}
Therefore,
\begin{equation*}
\frac{d\vec{a}}{dt} = (\nabla T_t^{-T}) \left(\frac{dn}{dt} - \frac{d(\nabla T_t^T)}{dt}\vec{a}\right).
\end{equation*}
Substituting back in the definition of $\vec{a}$ gives the alternative form of \eqref{eq:dTndt},
\begin{equation}\label{eq:dTndt_2}
\frac{d}{dt}(\nabla T_t^{-T}  n) = (\nabla T_t^{-T})\left( \frac{dn}{dt} - \frac{d (\nabla T_t^T)}{dt} (\nabla T_t^{-T}) n \right).
\end{equation}
Thie result in \eqref{eq:dTndt_2} is more advantageous than \eqref{eq:dTndt} because it only requires computing the 
derivative of $\nabla T_t^{T}$, whereas \eqref{eq:dTndt} requires computing the derivative of $\nabla T_t^{-T}$.

We compute the derivative of $\nabla T_t^T$ with respect to time $t$ next.  
%
Write $T_t(X) = T(t,X)$ to emphasize $T$ is a function of $X$ and $t$ 
and to allow for the subscript notation of partial derivatives.  
We adopt the notation $\nabla_x(\cdot)$ or $\nabla_X(\cdot)$ to denote the variable of spatial differentiation $x$ or $X$ when computing $\nabla T(t,\cdot)$.   
Computing the derivative with respect to time of $\nabla_X T(t,X)$ yields
\begin{align*}
\frac{d}{dt} \nabla_X T(t,X) &= \frac{d}{dt} \left(\partial_{X_j} T_i\right)_{1 \leq i, j \leq 3} 
	= \left(\partial_{X_j} \frac{dT_i}{dt}\right)_{1 \leq i, j \leq 3} 
	= \left(\partial_{X_j} V_i(X) \right)_{1 \leq i, j \leq 3}\\
	&= \nabla_X V(X). 
\end{align*}
%
Then,
\begin{align*}
\nabla_XV(X) = 
\begin{pmatrix} 
\DD\frac{\partial V_1}{\partial X_1} & \DD\frac{\partial V_1}{\partial X_2} & \DD\frac{\partial V_1}{\partial X_3} \\ 
\vdots & \ddots & \vdots \\
\DD\frac{\partial V_3}{\partial X_1} & \dots & \DD\frac{\partial V_3}{\partial X_3} 
\end{pmatrix}  
	&= \begin{pmatrix}
\DD\frac{\partial V_1}{\partial x_1} & \DD\frac{\partial V_1}{\partial x_2} & \DD\frac{\partial V_1}{\partial x_3} \\
\vdots & \ddots & \vdots \\
\DD\frac{\partial V_3}{\partial x_1} & \dots & \DD\frac{\partial V_3}{\partial x_3} 
\end{pmatrix}
\begin{pmatrix}
\DD\frac{\partial x_1}{\partial X_1} & \DD\frac{\partial x_1}{\partial X_2} & \DD\frac{\partial x_1}{\partial X_3} \\ 
\vdots & \ddots & \vdots \\
\DD\frac{\partial x_3}{\partial X_1} & \dots & \DD\frac{\partial x_3}{\partial X_3} 
\end{pmatrix} \\ 
	&= (\nabla_xV(x))(\nabla_XT(t,X)).
\end{align*}
This establishes 
\begin{equation}\label{eq:d_nablaT}
\frac{d}{dt}\nabla_X T(t,X) = (\nabla_x V(x))(\nabla_XT(t,X)).
\end{equation}
%
From \eqref{eq:d_nablaT} and the fact that the derivative of the transpose of a matrix is the transpose of the derivative,
\begin{equation}\label{eq:d_nablaT^T}
\frac{d}{dt}\nabla T_t^T = \left((\nabla V)(\nabla T_t)\right)^{T}.
\end{equation}
The equation \eqref{eq:d_nablaT^T} may be substituted in \eqref{eq:dTndt_2}, which may in turn be substituted in to \eqref{eq:dJ_s}, but we choose to substitute everything at the end of the section for neatness.

Next, we establish a formula for the normal derivative $\frac{dn}{dt}$ at time $t=0$.  Let $u$ and $v$ parameterize a surface in $\reals^3$ and define the surface by
\begin{equation}\label{eq:surface_R3}
\vec{r}(u,v) = (x(u,v), y(u,v), z(u,v)).
\end{equation}
With this parameterization of the surface, the transformation $T_t$ acting on the surface can be defined as
\begin{align}
\resizebox{.295\hsize}{!}{$T_t(x(u,v),y(u,v),z(u,v)) $}&= \resizebox{.65\hsize}{!}{$(x(u,v)+\alpha(u,v,t), y(u,v)+\beta(u,v,t), z(u,v)+\gamma(u,v,t))$} \notag\\
	&= \vec{r}(u,v) + \vec{s}(u,v,t), \label{eq:transformation_surface}
\end{align}
where $\vec{s}(u,v,t) = (\alpha(u,v,t),\beta(u,v,t),\gamma(u,v,t))$.  Note that $\vec{s}$ depends in $t$ but $\vec{r}$ is independent of $t$.  Next, define the (initial) normal vector to the surface in the reference coordinates by
\begin{equation}\label{eq:normal_surface}
n(X,0) = \DD\frac{r_u \times r_v}{|r_u \times r_v|},
\end{equation}
where the subscripts denote partial derivatives.  Then, at any time $t>0$, the transformed normal vector is calculated by applying \eqref{eq:transformation_surface} to the initial normal vector \eqref{eq:normal_surface},
\begin{equation}\label{eq:T_normal_surface}
n(X,t) = \DD\frac{(r_u+s_u) \times (r_v+s_v)}{|(r_u+s_u) \times (r_v+s_v)|}.
\end{equation}
For simplicity, let $\zeta = (r_u+s_u) \times (r_v+s_v)$.  Next, the initial time rate of change of \eqref{eq:T_normal_surface} 
is computed:
\begin{equation}  \label{eq:dkappa_dt}
\DD\left.\frac{dn(X,t)}{dt}\right|_{t=0} = \left.\DD\frac{d}{dt}\left(\DD\frac{\zeta}{|\zeta|}\right)\right|_{t=0} 
	= \left[\left.\frac{1}{|\zeta|}\frac{d\zeta}{dt} + \zeta\frac{d}{dt}\left(\frac{1}{|\zeta|}\right)\right]\right|_{t=0}. 
\end{equation}
Define 
\begin{equation*}
P = \frac{1}{|\zeta|}\frac{d\zeta}{dt} \qquad \textrm{and} \qquad Q = \zeta\frac{d}{dt}\left(\frac{1}{|\zeta|}\right).
\end{equation*}
The derivative $\DD\frac{d\zeta}{dt}$ can be computed using the product rule for cross products.
%
%
Notice that $\vec{r}$ does not depend on $t$, the simplified form of $P$ is
\begin{equation*}\begin{gathered}\begin{aligned}
\left.P\right|_{t=0} &= \left.\frac{1}{|\zeta|}\frac{d}{dt}((r_u+s_u) \times (r_v+s_v))\right|_{t=0}  \\
	&= \left.\frac{1}{|\zeta|}\left[\left(\frac{d}{dt}(r_u+s_u)\right) \times (r_v+s_v) + (r_u+s_u) \times \left(\frac{d}{dt}(r_v+s_v)\right)\right]\right|_{t=0} \\
	&= \left.\frac{1}{|(r_u+s_u) \times (r_v+s_v)|} \left[ \frac{ds_u}{dt} \times (r_v+s_v) + (r_u+s_u) \times \frac{ds_v}{dt} \right ] \right|_{t=0}\\
	&= \frac{1}{|r_u\times r_v|}\left [ \frac{ds_u}{dt}\times r_v + r_u \times \frac{ds_v}{dt} \right ]. 
\end{aligned}\end{gathered}\end{equation*}
The simplified form of $Q$ is 
\begin{align*}
\left. Q \right|_{t=0} &= \left. \zeta \frac{d}{dt}\left(\frac{1}{|\zeta|}\right) \right|_{t=0}
	= \left.\left((r_u+s_u) \times (r_v+s_v)\right)\frac{-1}{|\zeta|^3}(\zeta \cdot \frac{d\zeta}{dt}) \right|_{t=0}\\
	&= (r_u \times r_v) \frac{-1}{|r_u \times r_v|^3}\left[(r_u\times r_v) \cdot \left( \frac{ds_u}{dt}\times r_v + 
r_u \times \frac{ds_v}{dt} \right] \right).
\end{align*}
This gives the desired result:
\begin{align}
\frac{dn(X,0)}{dt} = & \frac{1}{|r_u\times r_v|}\left( \frac{ds_u}{dt}\times r_v + r_u \times \frac{ds_v}{dt} \right)  \nonumber \\
	& + (r_u \times r_v) \frac{-1}{|r_u \times r_v|^3}\left[(r_u\times r_v) \cdot \left( \frac{ds_u}{dt}\times r_v + 
r_u \times \frac{ds_v}{dt} \right) \right]. \label{eq:dTn_dt1}
\end{align}
To simplify, define 
\begin{equation*}
R =  \frac{ds_u}{dt}\times r_v + r_u \times \frac{ds_v}{dt}.
\end{equation*}
Then
\begin{equation*}\begin{gathered}\begin{aligned}
\frac{dn(X,0)}{dt} &= \frac{1}{|r_u\times r_v|} R  + (r_u \times r_v) \frac{-1}{|r_u \times r_v|^3}\left((r_u\times r_v) \cdot  R  \right) \\
	&= \frac{1}{|r_u\times r_v|} \left[R - \frac{r_u \times r_v}{|r_u \times r_v|} \left(\frac{r_u \times r_v}{|r_u\times r_v|} \cdot R\right)\right] \\
	&= \frac{1}{|r_u\times r_v|} \left[ R - n(X,0) \left(n(X,0)\cdot R\right)\right].
\end{aligned}\end{gathered}\end{equation*}
That is, 
\begin{equation}\label{eq:dTn_dt3}
\frac{dn(X,0)}{dt} = \frac{1}{|r_u\times r_v|} \left[ R - n(X,0) \left(n(X,0)\cdot R\right)\right].
\end{equation}
To finish the formula \eqref{eq:dJ_s} requires substituting \eqref{eq:d_nablaT^T} and an expression lengthier 
than \eqref{eq:dTn_dt1} into \eqref{eq:dTndt_2} and substituting the result into \eqref{eq:dJ_s}.  
The final expression for $\DD\frac{dJ_s}{dt}$ at arbitrary $t>0$ would be exceedingly long, and the expression is 
unnecessary for the calculations of the shape derivative. Indeed, only the derivative of $J_s$ at $t=0$ is used here.  
To compute this quantity, we notice that by \eqref{eq:d_nablaT^T} it follows 
\begin{equation}\label{eq:dT^T(0)}
\left.\frac{d}{dt}\nabla T_t^T\right|_{t=0} = ((\nabla V)(\nabla T_0))^T = (\nabla V)^T,
\end{equation}
where the $\nabla$ operator acting on $T_t^T(X)$ is with respect to $X$ and the one acting on $V(x)$ is with respect to $x$.  Now, using \eqref{eq:dJ_s}, \eqref{eq:dTndt_2}, \eqref{eq:dT^T(0)}, and \eqref{eq:dTn_dt3} gives the result
\begin{equation*}\begin{gathered}\begin{aligned}
\left.\frac{dJ_s}{dt}\right|_{t=0} &= (J_0(\nabla \cdot V)\circ T_0)\left|\nabla T_0^{-T} n\right| \\
	& \qquad + (\det\nabla T_0)\left(\frac{1}{\left|\nabla T_0^{-T}  n\right|}(\nabla T_0^{-T}  n)\cdot\left.\frac{d}{dt}(\nabla T_t^{-T}  n)\right|_{t=0} \right) \\
	&= \nabla \cdot V + \left(n\cdot \left.\frac{d}{dt}(\nabla T_t^{-T}  n)\right)\right|_{t=0} \\
	&= \nabla \cdot V + \left.\left( n \cdot \left[ (\nabla T_t^{-T}) \left( \frac{dn}{dt} - \frac{d (\nabla T_t^T)}{dt} (\nabla T_t^{-T}) n\right)  \right]\right) \right|_{t=0}  \\
	&= \nabla \cdot V + n \cdot \left( \frac{dn(X,0)}{dt} - (\nabla V)^T n\right) \\
	&= \nabla \cdot V + n \cdot \frac{1}{|r_u\times r_v|} \left( R - n \left(n\cdot R\right)\right)  - n \cdot (\nabla V)^Tn, 
\end{aligned}\end{gathered}\end{equation*}
where $n = n(X,0)$ is the initial unit normal, and hence $|n|=1$.  Notice that $n(n\cdot R)$ is the component of $R$ in the direction of $n$.  Then, $R-n(n\cdot R)$ is the component of $R$ perpendicular to $n$.  Thus, $n\cdot (R-n(n\cdot R)) = 0$.  This leads to 
\begin{equation}\label{eq:dJ_s(0)_1}
\left.\frac{dJ_s}{dt}\right|_{t=0} = \nabla \cdot V  - n \cdot (\nabla V)^Tn. 
\end{equation}
Next, notice that for any matrix $A$ and vector $\vec{x}$, $x \cdot A^Tx = x \cdot A x$.  This reduces \eqref{eq:dJ_s(0)_1} to
\begin{equation}\label{eq:dJ_s(0)_2}
\left.\frac{dJ_s}{dt}\right|_{t=0} = \nabla \cdot V  - n \cdot (\nabla V)n. 
\end{equation}
The right-hand side of \eqref{eq:dJ_s(0)_2} is exactly the surface divergence of $V$. This leads to the final expression 
of the initial time derivative of $J_s$,
\begin{equation}\label{eq:dJ_s(0)_b}
\left.\frac{dJ_s}{dt}\right|_{t=0} = \nabla_s \cdot V .
\end{equation}

\end{proof}



\end{document}